\documentclass[11pt,epsf]{article}
\usepackage{graphicx}
\usepackage{theorem, amsmath, amssymb}
\newtheorem{Theorem}{Theorem}[section]
\newtheorem{Definition}[Theorem]{Definition}
\newtheorem{Proposition}[Theorem]{Proposition}

\newtheorem{Lemma}[Theorem]{Lemma}

\theoremstyle{remark}
\newtheorem{Remark}[Theorem]{Remark}

  {\begin{equation}\left\{\begin{aligned}}%
  {\end{aligned}\right.\end{equation}\ignorespacesafterend}%

\newenvironment{SDE*}%
  {\begin{equation*}\left\{\begin{aligned}}%
  {\end{aligned}\right.\end{equation*}\ignorespacesafterend}%

\usepackage{color}
\begin{document}
\title{
 Stochastic Perron's method and elementary strategies  \\
for   zero-sum differential  games
}

\author{Mihai S\^{\i}rbu \footnote{University of Texas at Austin,
    Department of Mathematics, 1 University Station C1200, Austin, TX,
    78712.  E-mail address: sirbu@math.utexas.edu. The research of
    this author was supported in part by the National Science
    Foundation under Grant    DMS 1211988. Any opinions, findings, and conclusions or recommendations expressed in this material are those of the authors and do not necessarily reflect the views of the National Science Foundation.}}
\maketitle
\begin{abstract} We develop here the Stochastic Perron Method
in  the framework of  two-player  zero-sum  differential games.  We consider the formulation of the game where both players play, symmetrically, feed-back strategies
 (as  in   \cite{cr} and \cite{pz-game}) 
as opposed to the Elliott-Kalton formulation prevalent in the literature.
The class of feed-back strategies  we use is carefully chosen  so that 
   the state equation admits strong solutions
 and  the technicalities involved in  the Stochastic Perron Method carry through in a rather simple way. More precisely, we define the game  over  elementary strategies, which 
are well motivated by intuition.
Within this framework, the Stochastic Perron Method produces a viscosity sub-solution of the upper  Isaacs equation dominating the upper value of the game, and a viscosity super-solution of the upper Isaacs equation lying below the upper value of the game. 
Using a viscosity comparison result  we obtain that 
 the upper value is the unique and continuous viscosity solution of the upper Isaacs equation. An identical statement holds true for the lower value and the lower Isaacs equation.  A version of the Dynamic Programming Principle is obtained as a by-product.
  If the Isaacs condition is satisfied, the game has a  value over  elementary (pure) strategies. 
\end{abstract}


\section{Introduction}
Continuous time games have been an active research area since the pioneering work of Isaacs \cite{isaacs}.
 In \cite{ek}, Elliott and Kalton studied a deterministic differential game, where one player (using strategies)  observes continuously the actions of the other players, who plays open loop controls. This formulation of the game has been considered for the first time in stochastic framework in \cite{fs}. The viscosity solution approach previously used in deterministic setting turns out to be quite involved here, and, actually, one cannot prove the dynamic programming principle working  directly with the value functions (see \cite{fs} or \cite{MR1312585}). 

A different approach,  resembling verification for classical solutions but in the context of non-smooth viscosity solutions, was employed in \cite{swiech-1} and \cite{swiech-2} (first for deterministic problems and then for zero-sum stochastic games in Elliott-Kalton formulation). There, it was proven (using inf/sup convolution regularization of viscosity solutions in the spirit of \cite{MR791888} and an approximation of the state equation by non-degenerate diffusions) that  viscosity sub/super-solutions satisfy the corresponding half of the DPP (so called sub/super-optimality conditions), therefore comparing to the value function. The unique continuous viscosity solution (obtained  separately by PDE methods, most likely by   Perron's method for viscosity solutions introduced in  \cite{ishii}) is equal to the value function(s) and the DPP is satisfied. In order to regularize the diffusion, the probability space needs to be enlarged to accommodate an independent Brownian Motion. This means that verification is actually performed for a slightly different game, not the one in the original formulation.
Similar technical arguments to prove DPP for games in  Elliott-Kalton formulation were used very recently in \cite{krylov-1} and \cite{krylov-2} for non-degenerate diffusions, but with only measurable in space coefficients.



Out of the existing  literature on  games, most  papers (all above) use the Elliott-Kalton formulation  previously mentioned, where one player uses controls and the other uses strategies. This asymmetric formulation produces two values of the game that do not compare by definition, as one would expect. It seems more natural to consider a symmetric formulation of the game. One way to do this, is to consider, for both players, strategies which are based only on the knowledge of the past of the state. Each player's actions are responses to the other player's actions through the use of such  feed-back strategies. 
 Feed-back strategies have already been used in both one-player (control) problems (see the notion of natural controls in \cite{MR2723141}) and in (discretized) zero-sum deterministic games (under the name of positional/historical strategies, and used by one or both players) in \cite{krasovskii-subbotin-88}. From the recent literature on stochastic games, 
 \cite{cr}  consider a game in the strong formulation over feed-back strategies with delay, but even allow for some mixing, and  \cite{pz-game} use  discretized feed-back strategies (under the name of feed-back controls borrowed from the more classical control literature).
In the weak formulation of a non-Markovian game over feed-back strategies, \cite{pz-game} use BSDE's techniques to conclude that the values of the games are solutions to the path-dependent Isaacs equations. This work follows a long line of interesting results relating BSDEs to stochastic games  in Elliott-Kalton formulation:  \cite{hl}, \cite{ekh} or \cite{bl-1}, to name just a few.

We provide here a novel tool to study continuous-time zero-sum stochastic differential games, in an (arguably) more natural formulation. 
 The aim  of the current paper is two-fold:
\begin{enumerate}

\item to consider a symmetric formulation of the  stochastic  game over a (restricted class of) feed-back strategies that yield strong solutions to the state equations. The restriction is quite different, and (we believe) more natural than the delay in \cite{cr} or the deterministic time-discretization in \cite{krasovskii-subbotin-88} and \cite{pz-game}. We find it conceptually important to have a class of strategies where the state is a one-to-one response to the (unobserved) noise, once both players have decided on what strategies will be used,

\item to show that the Stochastic Perron Method (introduced 
 in \cite{bs-1} for linear problems, later  for Dynkin games  in \cite{bs-2} and for  control problems  in \cite{bs-3}) can be adapted to  differential games in the above formulation. This amounts to a new method to study games, and actually turns out to be rather elementary compared to the existing work.
\ 
\end{enumerate}
In order to illustrate our approach, we impose  some restrictive conditions on the controlled stochastic system, rather than consider the most general framework. More precisely: the state space is the whole $\mathbb{R}^d$, the state equation is  uniformly locally Lipschitz in the state variable, the terminal pay-off is  bounded, the controls live in some compact spaces, etc.
It is well known from \cite{fs} that stochastic games  are  non-trivial even with (more) restrictive assumptions, as the strategic nature of the problem remains intact. We also consider  only a terminal pay-off, but running cost can be easily added. However, we allow for  the state equation to  be degenerate, so one cannot expect any smoothness of solutions to the Bellman-Isaacs equation(s).

The method we present here  amounts to a verification result for non-smooth viscosity solutions.  The Stochastic Perron construction and viscosity comparison show that the value function(s) is the unique viscosity solution of the Bellman-Isaacs equation(s) and that  (a version of) the DPP holds.
Aside from the symmetric  definition of the game over feed-back strategies,
 our verification program is conceptually different  from the verification arguments in  \cite{swiech-1}, \cite{swiech-2}.
We work directly with the state equation and the original definition of the game. Our approach does  not  approximate the state equation by more regular diffusions so  {\em we have a  verification argument for the original (strong) definition of the game on  the initially chosen  physical space} that does not need to accommodate an additional independent Brownian motion.
  Technically, the Stochastic Perron method is rather elementary and is based only on how we define the stochastic semi-solutions, the well known  definition of viscosity solutions and It\^o's Lemma. Our verification is completed  by  a comparison argument that has to be proved analytically and separately anyway  (here and in other  existing approaches).

\section{Stochastic Differential Games and Strategies}
\subsection{The Stochastic System}
We consider a stochastic differential game with two players.
The first player's actions belong to a  compact metric space ($U, d_u)$
(usually $ U\subset \mathbb{R}^k$). The second player's actions belong also to a  compact metric space $(V,d_v)$ (usually   $V\subset \mathbb{R}^l$). We assume that the state lives in the whole space $\mathbb{R}^d$ (but one could also consider a natural  open domain $\mathcal{O}\subset \mathbb{R}^d$ as in \cite{bs-3}).
Let $b:[0,T] \times \mathbb{R}^d \times U \times V \to \mathbb{R}^d$ and $\sigma:[0,T]\times\mathbb{R}^d\times U \times V \to \mathbb{M}_{d,d'}$ be two continuous functions. We consider the diffusion 
\begin{equation}\label{eq:SDE}\left\{
\begin{array}{ll}
 dX_t=b(t,X_t,u_t, v_t)dt+\sigma (t, X_t,u_t, v_t)dW_t,   \\ 
 X_s=x \in \mathbb{R}^d,
 \end{array} \right.
 \end{equation}
 starting at an initial time $s$ at some  position $x$, and 
which is controlled by both players. Here, $W$ is a $d'$-dimensional Brownian motion on a \emph{fixed} probability space
 $(\Omega, \mathcal{F}, \mathbb{P})$. We intend to pose  the  problem  in such a way that the state equation will admit strong solutions. Therefore, is is natural to assume that, for the initial time $s$, the filtration considered is 
 $\mathbb{F}^s=(\mathcal{F}^s_t)_{s\leq t\leq T},$ the augmented natural filtration generated by the Brownian increments starting at $s$, by which we mean.
 $$\mathcal{F}^s_t=\sigma (W_u-W_s, s\leq u\leq t)\vee \mathcal{N}(\mathbb{P}, \mathcal{F})\ \  \textrm{for}\ \ s\leq t\leq T.$$
 However, any larger filtration that keeps $W_{\cdot}-W_s$ a Brownian motion starting at s and satisfying the usual conditions is also fine.
 Now, given a   bounded and continuous function $g:\mathbb{R}^d\rightarrow \mathbb{R}$, the second player pays to the first player the amount 
 $\mathbb{E}[g(X^{s,x;u,v}_T)].$ This leads to a zero-sum game of the form  
$$\sup _{u }\inf _v\mathbb{E}[g(X^{s,x;u,v}_T)],\ \ \ \inf _v\sup _u  \mathbb{E}[g(X^{s,x;u,v}_T)].\ \ \ \ \  \ 
$$
Just for the sake of presentation,  we have  assumed that no running cost is present, as in \cite{bs-3}. Such cost can easily be added. 
The lower and the upper Hamiltonians are defined by 

\[
H^-(t,x,p,M):=\sup_{u \in U} \inf _{v\in V}\left[b(t,x,u,v)\cdot p+\frac{1}{2}Tr \left (\sigma(t,x,u,v)\sigma(t,x,u,v)^T M\right)\right],\ \ 0\leq t\leq T, \; x\in \mathbb{R}^d,
\] 

\[
H^+(t,x,p,M):=\inf _{v\in V}\sup_{u \in U}\left[b(t,x,u,v)\cdot p+\frac{1}{2}Tr\left (\sigma(t,x,u,v)\sigma(t,x,u,v)^T M\right)\right],\ \ 0\leq t\leq T, \; x\in \mathbb{R}^d.
\] 
We associate, formally for the moment, the following Isaacs equations to the game
\begin{equation}\label{eq:Isaacs}
\left \{
\begin{array}{ll}
-v_t-H ^{i}(t,x,v_x,v_{xx})=0\ \ \textrm{on}\ [0,T)\times \mathbb{R}^d,\\
v(T,\cdot)=g(\cdot),\ \ \textrm{on}\ \mathbb{R}^d.
\end{array}
\right.
\end{equation}
For $i=-$ we have the lower Isaacs equation, and for $i=+$ we have the upper Isaacs equation.
We make the standing assumption that  that the coefficients $b, \sigma$ of the  stochastic system are jointly continuous, satisfy a uniform local Lipschitz condition in $x$, i.e.
\begin{equation*}\label{Lip}
{\bf (L)} \ \ |b(t,x,u,v)-b(t,y,u,v)|+|\sigma (t,x,u,v)-\sigma (t,y,u,v)|\leq L(K) |x-y|\ \ 
\end{equation*}
$\forall\  |x|,|y|\leq K, t\in [0,T],\  u\in U, v\in V$
for some $L(K)<\infty$, and 
also satisfy a global linear growth condition in $x$
$${\bf (LG)} \ \ \ 
|b(t,x,u,v)|+|\sigma(t,x,u,v)|\leq C(1+|x|)$$
$\forall\  |x|,|y|\in \mathbb{R}^d, t\in [0,T],\  u\in U, v\in V$ for some $C<\infty.$
Under these assumptions, if both players choose to hold a constant strategy, the controlled state has a unique solution. This holds true even if the game is started at a random (but stopping) time sequel to the starting time $s$, and the constant strategies depend on whatever happened before the starting time. More precisely, we have the following proposition:
\begin{Proposition} \label{prop:state-eq}
Let $s\leq \tau' \leq T$ be a stopping time with respect to  the filtration $\mathbb{F}^s=(\mathcal{F}^s_t)_{s\leq t\leq T}$. Let 
$\xi \in \mathbb{R}^d, a\in U, b \in V, $ be random variables measurable with respect to $\mathcal{F} ^s_{\tau '}.$ Then, the system starting at time $\tau '$ with initial condition $\xi$, where both players choose constant strategies $a,b$ in between $\tau '$ and $T$ has a (path-wise) unique strong solution. In other words, the SDE
\begin{equation}
\left \{
\begin{array}{ll}
dX_t=b(t,X_t, a, b)dt+\sigma (t, X_t,a, b)dW_t,\ \tau ' \leq t\leq T\\
X_{\tau '}=\xi
\end{array}
\right .
\end{equation}
has a unique  strong solution $(X_t)_{\tau ' \leq t\leq T}$. If $\xi$ is square integrable, the solution is square integrable as well.

\end{Proposition}
Proof: the arguments are standard, so we omit the details. The path-wise (or strong) uniqueness is proved by localization and Gronwall. Using again localization and then a Picard iteration, we prove the existence of a saturated strong solution, up to some explosion time. The global growth condition and Gronwall inequality insure that there is actually no explosion. $\diamond$
\subsection{Elementary Strategies}
So far, the game is defined formally. To have  a consistent mathematical definition, we need to define some  non-constant strategies for each player. In our  model both players observe continuously the state process (which is a response to both players actions) but not the noise, and 
both players play (symmetrically)  strategies based on the information available. 
This leads to a  formulation of the game  where the upper and the lower values compare by definition, as in in \cite{cr} and, more recently, in \cite{pz-game}. 
In other words,  we propose for use as strategies, by both players:

\begin{enumerate}
\item functionals of the past  of the state (pure feed-back strategies), if one expects the game to have a value. This is the case in \cite{cr} (with delay) and \cite{pz-game} (for the weak formulation),  following the older work  \cite{krasovskii-subbotin-88} on deterministic games, where feed-back strategies were discretized with respect to time,

\item strategies as above, but mixed (at each time), if the game over pure strategies as above does not have a value. Some mixing is actually allowed in \cite{cr} and  is also used  in \cite{krasovskii-subbotin-88} for the case of deterministic games.
\end{enumerate}
In other words, we believe that the game should be first formulated over such pure (feed-back) strategies and, if  the lower and upper values are not equal, then one should proceed directly  to considering mixed strategies where the game is expected to have a value.
For mathematical reasons, we treat in the present paper the  complete picture of the game over pure strategies. More precisely, we study separately the two value functions as  solutions of the two Isaacs equations. If the Isaacs condition is satisfied, then, obviously, the game is expected to have a value over pure strategies. We defer the  study  of elementary  mixed strategies  \emph{in the absence of  Isaacs condition} to forthcoming work \cite{sirbu-2}.


We do not assume any delay as in \cite{cr} for two reasons: such strategies don't work well with Stochastic Perron's Method and, conceptually, we find it harder to justify delaying the information available, while, at the same time, continuously changing actions. 
 Some other  natural restriction has still to be imposed on the feed-back strategies, so that, once each player chooses a strategy, the state equation has a unique  strong solution.  Despite the fact the noise is unobserved, we believe it is conceptually important for the state equation to have a unique strong solution over pure strategies. 
 In addition, our so called Elementary Strategies that we propose below, are exactly the kind of strategies that one can expect to enforce in a continuous-time game, where changing actions continuously may not be completely feasible, but when a change is made, the information available should not be delayed and only used later.  
 The  Elementary Strategies are defined in a very similar manner to how one defines  Elementary Integrands in stochastic integration, but on the state path space. 

 To simplify notation,  for the remainder of the paper, fixed a starting time $s$, we denote by $C([s,T])\triangleq C([s,T],\mathbb{R}^d)$ and endow this path space with the natural (and raw) filtration 
 $\mathbb{B}^s=(\mathcal{B}^s_t)_{s\leq t\leq T}$ defined by
 $$\mathcal{B}^s_t\triangleq \sigma (y(u),s\leq u\leq t), \ \ s\leq t\leq T.$$
 The elements of the path space $C([s,T])$ will be, usually, denoted by $y (\cdot)$ or $y$, when there is no risk of confusion. The stopping times on the space $C([s,T])$ with respect with the filtration $\mathbb{B}^s$, i.e.  mappings $\tau :C([s,T])\rightarrow [s,T]$ satisfying
 $  \{\tau \leq t\}\in \mathcal{B}^s_t \ \forall \ s\leq t\leq T$
are called  stopping  rules, following \cite{ks}. We denote   by $ \mathbb{B}^s$ the class of such stopping  rules starting at $s$.
 \begin{Definition}[Elementary Pure Strategies] \label{def:s} Fix $0\leq s\leq T$.
 An  elementary  strategy $u$ starting at $s$,  for the first player,  is defined by

\begin{itemize}
\item a finite non-decreasing sequence of stopping  rules, i.e. $\tau _k \in \mathbb{B}^s$ for $k=1, \dots, n$ and 
$$s=\tau _0\leq \dots \tau _k\leq \dots \leq \tau _n=T $$
\item for each $k=  1\dots n$, a constant value of the strategy $\xi_k$ in between the times $\tau _{k-1}$ and $\tau _k$, which is decided based only on the knowledge of the past state up to  $\tau _{k-1}$, i.e.
$\xi_k:C([s,T])\rightarrow U$ such that
$\xi _k\in \mathcal{B}^s_{\tau _{k-1}}$.
\end{itemize}
 The strategy is to hold $\xi_k$ in between $(\tau _{k-1}, \tau _{k}]$, i.e.
 $u: (s,T]\times C([s,T])\rightarrow U$ is defined by
 $$u(t, y(\cdot))\triangleq \sum _{k=1}^n \xi _{k}(y(\cdot))1_{\{ \tau _{k-1}(y(\cdot))<t\leq  \tau _{k}(y(\cdot))\}}.$$
Such a strategy can be denoted, when there is no risk of confusion, also as
$u=(\xi _k, \tau _k)_{k=1,n}$.
  An elementary strategy $v$ for the second player is defined in an identical way, but takes values in $V$.
 We denote by $\mathcal{U}(s,s)$ and $\mathcal{V}(s,s)$ the collections of all possible elementary strategies for the first, and the second player, given the initial deterministic time $s$.
\end{Definition}
The timings $\tau _k$ of changing the controls in the definition of strategies, are exactly the  "stopping rules"  considered in the game of stopping and control in \cite{ks}.  In addition, we would like to point out that, since $\xi_1\in \mathcal{B}^s_s$, this means that 
$\xi_1(\cdot)=l(y(s)),$ for some measurable function $l:\mathbb{R}^d\rightarrow U.$
It is obvious that, for a simple strategy, the representation $u=(\xi _k, \tau _k)_{k=1,n}$ is not unique.
\begin{Remark}The elementary strategies  are obviously non-anticipative functionals of the paths of $X$ but  they are neither
\begin{enumerate} 
\item Lipschitz with respect to the sup norm on the path-space  (this allows for great flexibility in pursuing Stochastic Perron)
\item nor delayed strategies as considered in \cite{cr}
\end{enumerate} 
Lipschitz functionals of the path are good strategies to define the values of the game (resulting in strong solutions of the state equation, with the very strong  additional Lipschitz assumption of coefficients) but  hard to fit in the framework of Stochastic Perron we develop here. In addition, we believe that our  "discrete-time" strategies are easier to justify in applications. 
\end{Remark}
Fortunately, the strategies we just defined still produce strong solutions of the state equation without any  Lipschitz assumption in $u,v$. More precisely, iterating Proposition \ref{prop:state-eq}, we get
\begin{Proposition}\label{prop:state-eq-simple} Fix $s,x$ and 
  let  players one and two choose strategies $u\in \mathcal{U}(s,s)$ and $ v \in\mathcal{V}(s,s)$. Then, there exists a unique strong (and square integrable) solution
$(X^{s,x;u,v}_t)_{s \leq t\leq T}$ of the state equation
\begin{equation}
\left \{
\begin{array}{ll}
dX_t=b(t,X_t,  u(t, X_{\cdot}), v(t, X_{\cdot})dt+\sigma (t, X_t,u(t, X_{\cdot}), v(t, X_{\cdot}))\, dW_t,\ s \leq t\leq T\\
X_{s}=x\in \mathbb{R}^d.
\end{array}
\right .
\end{equation}
Uniqueness holds in the path-wise sense.
\end{Proposition}
The proof of the above proposition is based on iterating Proposition \ref{prop:state-eq} together with the following very simple but useful lemma:
\begin{Lemma}\label{lemma:stopping} Fix $s$ and 
let $\tau$ be a stopping  rule,
$\tau :C([s,T])\rightarrow [s,T]$, $\tau \in \mathbb{B}^s$. Let $(X_t)_{s\leq t\leq T}$  be a process with continuous (all, not only almost surely) paths, which is adapted to $\mathbb{F}^s.$ Then,  the random time
$\tau _X:\Omega \rightarrow [s,T]$ defined by
$\tau _X(\omega)\triangleq \tau (X_{\cdot}(\omega))$ is a stopping time w.r.to the filtration $\mathbb{F}^s.$ In addition
$X_{\tau _X}\in \mathcal{F} ^s_{\tau _X}.$
\end{Lemma}
\begin{Remark} In one-player (control) problems, the  use of open-loop controls yields strong solutions of the state equation (because of the  Lipschitz condition in the standing assumptions). Something similar would happen if one defined the game over open-loop controls, but this formulation is not appropriate for games.  In the Elliott-Kalton formulation of \cite{fs}, where one player chooses open-loop controls and the other strategies, again, existence of strong solutions of the state equation is easily checked. Once feed-back strategies are allowed (actually for both players), the existence of strong solutions is non-trivial. One can, therefore, either restrict strategies, as in \cite{cr} in order to have a strong formulation, or consider a weak formulation as in \cite{pz-game} (but restriction are needed even in the weak formulation). We choose here  the strong formulation, and restrict the strategies to the class of elementary strategies. The restriction to the class of elementary strategies is meaningful both at the modeling level (actions can only be changed discretely in time) and at the level of the mathematics (we need the state equation to be well posed).
\end{Remark}
We  now formulate the game rigorously. For  fixed $0\leq s \leq t \leq T$ and $x\in \mathbb{R}^d$ (deterministic), we define the lower and the upper value of the game
$$V^-(s,x)\triangleq \sup _{u \in  \mathcal{U} (s,s) }\inf _{v\in \mathcal{V} (s,s)} \mathbb{E}[g(X^{s,x:u,v}_T)]\leq 
\inf _{v\in \mathcal{V}(s,s)}\sup _{u\in \mathcal{U}(s,s)}  \mathbb{E}[g(X^{s,x;u,v}_T)]\triangleq V^+(s,x).$$
  Since both our players play strategies (symmetrically), the lower and the upper values compare \emph{by definition}. This feature has already been emphasized in the literature, like \cite{cr} or {\cite{pz-game}, in opposition  to  the Elliott-Kalton formulation in \cite{fs}. 
 A finite valued version of  feed-back strategies is present in \cite{pz-game} for technical reasons, but the discretization  is over deterministic times, and the game is considered in weak formulation. Our strategies are actually not finitely  valued: they are discrete  only in time. The discretization over deterministic times in \cite{pz-game} is also reminiscent of the classic work \cite{MR945913}, where the author uses a fixed time grid and then passes to the limit using semi-group techniques.
  
  The lower value of the game has the meaning that player 1 announces his/her full strategy to player 2, who chooses a strategy accordingly. Symmetrically for the upper value. Of course, something like this does not really happen, so, if the two values differ, one should  just consider mixed strategies. However, we do analyze the two possibly different values for mathematical reasons. Assume for a moment, that, player 2, indeed, knows the full strategy that player 1 is going to  use. If that is the case, since both players' strategies are non-anticipative functionals of the path, conceptually, this would correspond exactly to the Elliott-Kalton idea that player 2 observes continuously player's 1 controls. It is, therefore, not surprising that we expect such lower value and upper value, to be the solutions to the corresponding Isaacs equation. However, such a heuristic statement needs a proof, and part of our contribution is to provide such proof.

We also need  need to define strategies starting at  sequel  times to the initial (deterministic) time $s$. The starting time is  a stopping rule.
\begin{Definition}[Elementary Strategies starting later] \label{def:tau}Fix $s$ and let 
 $\tau \in \mathbb{B}^s$ be  a stopping  rule.  An  elementary strategy, denoted by  $u\in \mathcal{U}(s,\tau)$,  for the first player, starting at $\tau$, is defined by

\begin{itemize}
\item (again) a finite non-decreasing sequence of stopping  rules, i.e.
$\tau _k \in \mathbb{B}^s,$ $k=1,\dots n$ 
for some finite $n$, and with 
$\tau =\tau _0\leq \dots \tau _k\leq \dots \leq \tau _n=T.$
\item for each $k=  1\dots n$, a constant action $\xi_k$ in between the times $\tau _{k-1}$ and $\tau _k$, which is decided based only on the knowledge of the past state up  $\tau _{k-1}$, i.e.
$\xi_k:C([s,T])\rightarrow U$ such that
$\xi _k\in \mathcal{B}^s_{\tau _{k-1}}$.
\end{itemize}
 The strategy is, again, to hold $\xi_k$ in between $(\tau _{k-1}, \tau _{k}]$, i.e.. 
 $$u: \{(t,y)| \tau (y)<t\leq T, y\in  C([s,T])\}\rightarrow U \ \textrm{
 with }\ \ 
 u(t, y(\cdot))\triangleq \sum _{k=1}^n \xi _{k}(y(\cdot))1_{\{ \tau _{k-1}(y(\cdot))<t\leq  \tau _{k}(y(\cdot))\}}.$$
 
 We define similarly $\mathcal{V}(s,\tau)$ for the second player.
 \end{Definition}
 It is clear that, if we consider $\tau=s$   in the Definition \ref{def:tau}  of $\mathcal{U}(s, \tau)$, we obtain the same thing as  $\mathcal{U}(s,s)$ in Definition \ref{def:s}, so the notation is consistent. The same holds true for $\mathcal{V}$. The classes of strategies $\mathcal{U}(s,\tau)$ and $\mathcal{V}(s,\tau)$ are, obviously, not closed in any reasonable sense, but they are  robust.  The two results below have very simple proofs, which we omit.
\begin{Lemma} Let $u_1, u_2\in \mathcal{U}(s,\tau)$ for some $s\leq T$ and $\tau \in \mathbb{B}^s$. Let $A\subset C([s,T]$ such that $A\in \mathbb{B}^s_{\tau}$. Then, the mapping 
 $u :(s,T]\times C([s,T])\rightarrow U $ defined by
$$u(t, y)=1_{\{y\in A\}}u_1(t, y)+1_{\{y\in A^c\}}u_2 (t, y),$$
is an elementary strategy in $\mathcal{U}(s,\tau)$. An identical statement holds for the second player.
\end{Lemma}
 Strategies in $\mathcal{U}(s,\tau)$ or $\mathcal{V}(s, \tau)$ cannot be used by themselves for the game starting at $s$, but have to be concatenated with other strategies. 
 \begin{Proposition}[Concatenated elementary strategies]\label{prop:conc} Fix $s$ and 
 let $\tau \in \mathbb{B}^s$ be  a stopping  rule  and $\tilde{u}\in \mathcal{U}(s,\tau).$ Then, for each $u\in \mathcal{U} (s,s)$, the mapping
 $u \otimes_{\tau} \tilde{u}:(s,T]\times C([s,T])\rightarrow U $
 defined by 
 $$\big (u \otimes _{\tau} \tilde{u} \big )(t,y(\cdot))\triangleq u(t, y(\cdot)) 1_{\{s<t\leq \tau (y(\cdot))\}}+\tilde{u}(t, y(\cdot))1_{\{\tau (y(\cdot))<t\leq T\}}$$
 is a simple strategy starting at $s$, i.e. 
 $u \otimes _{\tau} \tilde{u}\in \mathcal{U}(s,s).$
 A similar statement holds for the second player.  \end{Proposition}

\section{Stochastic Perron's Method}
Having defined elementary strategies, and the concatenation,  we are now in position to define stochastic semi-solutions of the two  Isaacs equations.
\begin{Definition}[Stochastic Super-Solution of Upper Isaacs]
A function $w:[0,T]\times \mathbb{R}^d \rightarrow \mathbb{R}$ is called a stochastic super-solution of the upper Isaacs equation if
\begin{enumerate}
\item it is bounded, continuous  and $w(T, \cdot)\geq g(\cdot)$,
\item  for each $s$ and for each stopping  rule $\tau :C([s,T])\rightarrow [s,T]$, $\tau \in \mathbb{B}^s$
 there exists   an elementary strategy $\tilde {v}\in \mathcal{V}(s,\tau)$ such that, for any
$u\in \mathcal{U}(s,s)$,  any $v\in \mathcal{V}(s,s)$, any $x\in \mathbb{R}^d$ and  each stopping  rule  $\rho \in \mathbb{B}^s$, $\tau\leq \rho\leq T$,
with the simplifying  notation
$X\triangleq X^{s,x,u, v \otimes _{\tau} \tilde v}$ and (see Lemma \ref{lemma:stopping})
$\tau '\triangleq \tau (X), \rho '\triangleq \rho (X),$
we have 
$$w(\tau ', X_{\tau '})\geq \mathbb{E}[w (\rho', X_{\rho'})  |\mathcal{F}^s_{\tau'}]\ \ \mathbb{P}- a.s.$$
\end{enumerate}
\end{Definition}
Choosing $\tau =s$, we can see that, if $w$ is a stochastic super-solution, there exists 
$\tilde{v}\in \mathcal{V}(s,s)$ such that
$$w(s,x)\geq \mathbb{E} \left [w  (\rho (X^{s,x,u, \tilde{v}}_{\cdot}), X^{s,x, u,  \tilde{v}}_{\rho (X^{s,x,u, \tilde{v}}_{\cdot})} )|\mathcal{F}^s_s \right], \ \mathbb{P}-a.s$$
for all $u\in \mathcal{U}(s,s)$ and  $\rho \in \mathbb{B}^s$.
After taking the expectation, it is now obvious that, if $w$ is a stochastic super-solution, then  we have the  half DPP for $w$, i.e.
\begin{equation}\label{dpp-upp-upp}w(s,x)\geq \inf _{v\in \mathcal{V}(s,s)} \sup _{u\in \mathcal{U}(s,s)}\mathbb{E}\left [w(\rho(X^{s,x,u,v}_{\cdot}), X^{s,x,u,v}_{ \rho (X^{s,x,u,v}_{\cdot}   )})\right ], \ \ \ \forall \rho \in \mathbb{B}^s.
\end{equation} Since $w(T, \cdot)\geq g(\cdot)$, we obtain easily that
$w(s,x)\geq V^+(s,x).$
Compared to the Elliott-Kalton formulation prevalent in the literature, we have here a game where both players use elementary  feed-back  strategies. Also, in \eqref{dpp-upp-upp}, stopping rules are considered instead of stopping times. Therefore,  even if one tried to enlarge the probability space to accommodate an independent Brownian motion, it is not clear if the non-trivial approximation arguments in \cite{swiech-2} could be reproduced to conclude that the super-optimality principle \eqref{dpp-upp-upp} is equivalent to the notion of viscosity super-solution of upper Isaacs equation.
However, instead of following this direction,  we circumvent  most difficulties using the more elementary Stochastic Perron Method, which means  taking the infimum of stochastic super-solutions.

The stochastic sub-solutions of the lower Isaacs equation are defined symmetrically, in an obvious way.
\begin{Definition}[Stochastic Sub-Solution of Lower Isaacs]
A function $w:[0,T]\times \mathbb{R}^d \rightarrow \mathbb{R}$ is called a stochastic sub-solution of the lower  Isaacs equation if
\begin{enumerate}
\item it is bounded, continuous  and $w(T, \cdot)\leq g(\cdot)$,
\item for each $s$ and  for each stopping  rule $\tau \in \mathbb{B}^s$ 
 there exists   an elementary strategy $\tilde {u}\in \mathcal{U}(s,\tau)$ such that, for any
$u\in \mathcal{U}(s,s)$,  any $v\in \mathcal{V}(s,s)$, any $x$ and  each stopping  rule  $\rho \in \mathbb{B}^s$, $\tau\leq \rho\leq T$,
with the simplifying  notation
$X\triangleq X^{s, x,u \otimes _{\tau} \tilde{u}, v }$ and 
$\tau '\triangleq \tau (X), \rho '\triangleq \rho (X),$
we have 
$$w(\tau ', X_{\tau '})\leq \mathbb{E}[w(\rho', X_{\rho'})|\mathcal{F}^s_{\tau'}]\ \ \mathbb{P}- a.s.$$
\end{enumerate}\end{Definition}
It is another easy observation that, a  stochastic sub-solution of the lower  Isaacs equation satisfies
$w\leq V^-$ and the corresponding half DPP/sub-optimality principle.  The two definitions above are perfectly symmetric, and they would be enough to proceed with Stochastic Perron, in case the Isaacs condition holds. For the general case we  study here we need to also define:

\begin{Definition}[Stochastic Sub-Solution of Upper Isaacs]
A function $w:[0,T]\times \mathbb{R}^d \rightarrow \mathbb{R}$ is called a stochastic sub-solution of the upper  Isaacs equation if
\begin{enumerate}
\item it is bounded, continuous  and $w(T, \cdot)\leq g(\cdot)$,
\item for each $s$ and  for each stopping rule $\tau :C([s,T]\rightarrow [s,T]$, $\tau \in \mathbb{B}^s$ and each strategy $v\in \mathcal{V}(s,s)$
 there exists   an elementary  strategy $\tilde {u}\in \mathcal{U}(s,\tau)$ (depending on $v$ and $\tau$) such that, for any 
$u\in \mathcal{U}(s,s)$  and  any $x$ as well as each stopping rule   $\rho \in \mathbb{B}^s$, $\tau\leq \rho\leq T$, with the simplifying  notation
$X\triangleq X^{s, x,u \otimes _{\tau} \tilde{u}, v }$ and 
$\tau '\triangleq \tau (X), \rho '\triangleq \rho (X),$
we have 
$$w(\tau ', X_{\tau '})\leq \mathbb{E}[w(\rho', X_{\rho'})|\mathcal{F}^s_{\tau'}]\ \ \mathbb{P}- a.s.$$
\end{enumerate}\end{Definition}
We omit the definition of  Stochastic super-solution of lower Isaacs for the reason of being obvious by now. 
Let $w$ a stochastic sub-solution of upper Isaacs. Fix $v\in \mathcal{V}(s,s)$ and $\tau =s$. There exists $\tilde{U}\in \mathcal{U}(s,s)$ (depending on $v$) such that,  for each $\rho \in \mathbb{B}^s$ we have 
$$w(s,x)\leq \mathbb{E}\left [w(\rho (X^{s,x,\tilde{u}, v}_{\cdot}), X^{s,x, \tilde{u},  v}_{\rho (X^{s,x,\tilde{u}, v}_{\cdot})})|\mathcal{F}^s_s   \right], \ \mathbb{P}-a.s$$
After taking the expectation, it is now obvious that, if $w$ is a stochastic super-solution, then  we have the half DPP/sub-optimality principle for $w$, i.e.
\begin{equation}\label{dpp-upp-low}
w(s,x)\leq \inf _{v\in \mathcal{V}(s,s)} \sup _{u\in \mathcal{U}(s,s)}\mathbb{E}\left [w(\rho(X^{s,x,u,v}_{\cdot}), X^{s,x,u,v}_{ \rho (X^{s,x,u,v}_{\cdot}   )})\right ], \ \ \ \forall \rho \in \mathbb{B}^s.
\end{equation} Since $w(T, \cdot)\leq g(\cdot)$, we obtain easily that
$w(s,x)\leq V^+(s,x).$
\begin{Remark} The definitions of stochastic semi-solutions are quite different from those in the case of one player (control problems) in \cite{bs-3}. There, one simply starts the problem at a stopping time sequel to $s$ in the definition(s): the same does not make sense here. Therefore, the   nearly optimal strategy  (maybe depending on the other player strategy, as  for the sub-solution of upper Isaacs) in the definition(s) has to be concatenated with  any other strategy up to that stopping time rule $\tau$. This makes sense conceptually as one cannot simply start the game at $\tau$. On the other hand, if we tried to start the game at a stopping time $\tau '$ of the original filtration $\mathbb{F}^s$, we would get into significant difficulty trying  to complete Stochastic Perron's method below. In some sense, this is one of the best features of Stochastic Perron's Method, as it allows for great flexibility as to what we can choose for stochastic semi-solutions, as long as the definitions provide an easy comparison to the value function(s).
\end{Remark}


We denote by

1.  $\mathcal{U}^+$ the set of stochastic super-solutions of upper Isaacs

2.  $\mathcal{U}^-$ the set of stochastic  sub-solutions of upper Isaacs

3. $\mathcal{L}^+$ the set of stochastic super-solutions of lower Isaacs (not defined explicitly)

4. $\mathcal{L}^-$ the set of stochastic sub-solutions of lower Isaacs.

Since  $g$ is bounded, all the  sets are obviously non-empty. As expected, by Stochastic Perron's Method we mean taking the sup of sub-solutions and/or inf of super-solutions (by which we mean stochastic semi-solutions):
$$v^- \triangleq \sup _{w\in \mathcal{U}^-} w\leq  V^+\leq \inf _{w\in \mathcal{U}^+}w\triangleq v^+,$$
and
$$w^- \triangleq \sup _{w\in \mathcal{L}^-} w\leq  V^-\leq \inf _{w\in \mathcal{L}^+}w\triangleq w^+ .$$
We can say, without  need for any more proof, that
$w^-\leq V^-\leq V^+\leq v^+.$ We recall the standing assumptions: $g$ is continuous and bounded, $b$ and  $\sigma$ are continuous on their whole corresponding domains, locally  uniformly Lipschitz in $x$ and have global linear growth in $x$.
\begin{Theorem}[Stochastic Perron for games]\label{main1} Under the standing assumptions:
\begin{enumerate}
\item 
the function $v^+$ is a bounded upper semi-continuous (USC) viscosity sub-solution of the  upper Isaacs equation and satisfies the Half DPP \eqref{dpp-upp-upp}. The function $v^-$ is a bounded lower semi-continuous (LSC) viscosity super-solution of the  upper Isaacs equation and satisfies the half DPP \eqref{dpp-upp-low}.
\item the function $w^+$ is a bounded upper semi-continuous (USC) viscosity sub-solution of the  lower Isaacs equation, and the function $w^-$ is a bounded lower semi-continuous (LSC) viscosity super-solution of the  lower Isaacs equation (and they both satisfy the corresponding halves of the DPP for the lower equation).
\end{enumerate}\end{Theorem}
\begin{Remark} The viscosity semi-solution property holds at  time $t=0$ as well, consistent with the notation in \eqref{eq:Isaacs}. This is important (as can be seen from Lemma \ref{lem:visc}) to obtain comparison up to  time $t=0$ for semi-continuous viscosity solutions.
\end{Remark}
In order to prove the Theorem, we need  two Lemmas.
\begin{Lemma}
\begin{enumerate}
\item if $w_1, w_2\in \mathcal{U}^+$ then $w_1\wedge w_2\in \mathcal{U}^+$
\item if $w_1,w_2\in \mathcal{U}^-$ then $w_1\vee w_2\in \mathcal{U}^-$
\end{enumerate}
\end{Lemma}
Fix $\tau \in \mathbb{B}^s$  a stopping time strategy and. Let $\tilde{v}_1, \tilde{v}_2 \in \mathcal{V}(s, \tau)$  be the two strategies, for the second player, starting at $\tau$ corresponding the the super-solutions $w_1$ and $w_2$. The new strategy starting at $\tau$ defined by 
$$
\tilde{v}(t, y(\cdot))=\tilde{v}_1 (t, y (\cdot)) \, 1_{\{w_1(\tau (y), y (\tau (y)))\leq w_2(\tau (y), y (\tau (y)))\}}+\tilde{v}_2  (t, y(\cdot))\, 1_{\{w_1(\tau (y), y (\tau (y)))> w_2(\tau (y), y (\tau (y)))\}}$$
does the job for the definition of $w\triangleq w_1\wedge w_2$ as a stochastic super-solution of the upper Isaacs.  A similar  construction has to be made for the second part.  $\diamond$
\begin{Lemma}There exists a non-increasing sequence $\mathcal{U}^+\ni w_n\searrow v^+$ and an non-decreasing sequence $\mathcal{U}^-\ni v_n\nearrow v^-.$
\end{Lemma}
Proof: according to Proposition 4.1 in \cite{bs-1}, there exist $\tilde{w}_n\in \mathcal{U}^+$ such that
$v^+=\inf _n \tilde{w}_n.$ Now, we can just define
$w_n =\tilde{w}_1\wedge \dots \wedge \tilde{w}_n\in \mathcal{U}^+\searrow v^+.$ A very similar proof works for the second part of the lemma. $\diamond$

\noindent Proof of Theorem \ref{main1}:  the proof of the half DPP is obvious, by  construction, since we have seen that the stochastic semi-solution satisfy such half DPP. As far as the viscosity property is concerned,  we only prove item 1, as the second is obviously symmetric. 


\noindent {\bf 1.\ \ ( $v^+$ is viscosity sub-solution of upper Isaacs) }

1.1:  The interior sub-solution property for $v^+$: Let $(t_0,x_0)$ in the parabolic  interior $[0,T)\times \mathbb{R}^d$ such that a smooth function $\varphi$ strictly touches $v^+$ from above at $(t_0,x_0)$. Assume, by contradiction, that
$$\varphi _t +H^+(t,x,\varphi _x, \varphi _{xx})<0\ \  \textrm{ at}\ \  (t_0, x_0).$$ 
In particular, there exists $\hat{v}\in V$ and $\varepsilon >0$ such that
$$\varphi _t (t_0,x_0)+ \sup_{u \in U}\left[b(t_0, x_0,u, \hat{v})\cdot \varphi _x(t_0,x_0)+\frac{1}{2}Tr(\sigma(t_0,x_0,u,\hat{v})\sigma(t,x,u,\hat{v})^T \varphi _{xx}(t_0,x_0)) \right] <- \varepsilon.$$
To simplify notation, all small balls here are actually included in (i.e. intersected with)   the parabolic interior.
 Since $b, \sigma$  are continuous, and $U$ is compact,  the uniform continuity of the above expression in $(t,x,u)$ for $(t,x)$ around $(t_0, x_0)$ implies that  there exists a smaller $\varepsilon >0$ such that
$$\varphi _t(t,x)+ \sup_{u \in U}\left[b(t, x,u, \hat{v})\cdot \varphi _x(t,x)+\frac{1}{2}Tr(\sigma(t,x,u,\hat{v})\sigma(t,x,u,\hat{v})^T \varphi _{xx}(t,x)) \right] <-\varepsilon,\ \ {\textrm on}\ B(t_0, x_0, \varepsilon).$$
Now, on the compact (rectangular) torus $\mathbb{T}= \overline{B(t_0, x_0, \varepsilon)}- B(t_0, x_0, \varepsilon/2)$ we have that $\varphi >v^+$ and the min of $\varphi -v^+$ is attained, therefore it is strictly positive. In other words
$\varphi >v^+ +\eta$ on $\mathbb{T}$ for some $\eta >0$. Since $w_n\searrow v^+$, a Dini type argument similar to \cite{bs-2} and \cite{bs-3} shows that, for $n$ large enough we have $\varphi >w_n+\eta/2$. For simplicity, fix such an $n$ and call $w=w_n$. Now, define, for small $\delta <<\eta /2$
$$w^{\delta}\triangleq
\left \{ 
\begin{array}{ll}
(\varphi -\delta)\wedge w\ \ {\textrm on}\ \  B(t_0, x_0, \varepsilon),\\
w \ \ \textrm{outside}\ \  B(t_0, x_0, \varepsilon).
\end{array}
\right.
$$
Since $w^{\delta}(t_0,x_0)<v^+(t_0,x_0)$, we obtain a contradiction if we can show $w^{\delta}\in \mathcal{U}^+$. 

In order to do so, fix $s$ and let $\tau\in \mathbb{B}^s $ be  a stopping rule for the  initial time $s$. We need to construct an elementary  strategy $\tilde{v}\in \mathcal{V}(s,\tau)$ in the Definition of stochastic super-solution for $w^{\delta}$. This can be done  as follows:
since $w$ is a stochastic super-solution of upper Isaacs, there exists  an elementary strategy $\tilde{v}_1$ for $w$ starting at $\tau$ that does the job. Then
\begin{enumerate}
\item if $(\varphi -\delta) <w$ at $\tau$, follow the constant action $\hat{v}.$
\item if $(\varphi -\delta) \geq w$ at $(\tau, \xi)$ follow the strategy $\tilde{v}_1$
\item follow the strategy  defined in 1-2 until the first time $\tau _1$ when $(t,X_t)\in \partial B(t_0, x_0, \varepsilon/2)$. On this boundary,  we know that $w^{\delta}=w$.
\item after this, follow the strategy $\tilde{v}_3\in \mathcal{V}(s, \tau _1)$ corresponding to the stochastic sub-solution $w$ with starting stopping rule $\tau _1$
\end{enumerate}
In a rigorous notation though, we first define $\tilde{v}_2\in \mathcal{V}(s,\tau)$ by 
$$\tilde{v}_2(t,y (\cdot))= \hat{v} \, 1_{\{\varphi(\tau (y), y (\tau (y))) -\delta < w(\tau (y), y (\tau (y)))\}}+\tilde{v}_1 (t, y(\cdot))\, 1_{\{\varphi(\tau (y), y (\tau (y))) -\delta  \geq w(\tau (y), y (\tau (y)))\}}.$$
Denote now by 
$\tau _1:C([s,T])\rightarrow [s,T]$, $\tau \leq \tau _1\leq T$, 
$$\tau _1(y)\triangleq\inf _{\tau (y)\leq t\leq T}\{(t,y(t))\in \partial B(t_0, x_0, \varepsilon/2)\}.$$
Now, if $\tilde{v}_3\in \mathcal{V}(s,\tau _1)$ is the strategy from the definition of  the stochastic super-solution $w$, corresponding to stopping rule $\tau _1$, then we can define
$$\tilde{v}=\tilde{v}_2\otimes _{\tau _1}\tilde{v}_3\in \mathcal{V}(s, \tau)$$
to finish the proof. In order to check that $\tilde{v}$ fulfills the condition in the definition of stochastic super-solution for upper Isaacs, one has to follow similar arguments to the proofs in \cite{bs-3}, based on It\^o's Lemma applied to $d(\varphi-\delta)(t, X_t)$, conditioning and the tower property. However, one has to shift attention from stopping times as in \cite{bs-3} to stopping rules. We present the arguments in what follows. Fix  
$u\in \mathcal{U}(s,s)$ , $v\in \mathcal{V}(s,s)$,   $x\in \mathbb{R}^d$ and  choose a stopping  rule   $\rho \in \mathbb{B}^s$ with $\tau\leq \rho\leq T$. Denote by  
$X\triangleq X^{s,x,u, v \otimes _{\tau} \tilde v}$ , where $\tilde{v}$ was just defined above and $\tau '\triangleq \tau (X), \rho '\triangleq \rho (X).$  Let also 
$\tau '_1\triangleq \tau _1(X)$ such that
$\tau '\leq \tau '_1\leq \rho '$, and define the event
$A\triangleq \{(\varphi -\delta )(\tau ', X_{\tau '})< w(\tau ', X_{\tau '})\}$. We first observe that
$$X_t=1_A  X_t^{s,x,u, v \otimes _{\tau} \hat{ v}}+1_{A^c} X_t^{s,x,u, v \otimes _{\tau} \tilde v_1}, \tau '\leq t\leq \tau '_1.$$
On $A$, the process $(\varphi -\delta)(\cdot, X_{\cdot})$ then satisfies the super-martingale property in between $\tau '$ and $\tau '_1$ from It\^o's Lemma. On $A^c$, the process  $w(\cdot, X_{\cdot})$ satisfies the super-martingale inequality  in between $\tau '$ and $\tau '_1$ from the definition of the strategy $\tilde{v}_1\in \mathcal{V}(s, \tau).$ Putting this together, we have 
$$w^{\delta}(\tau ', X_{\tau '})\geq \mathbb{E}[w^{\delta}(\tau _1', X_{\tau _1'})|\mathcal{F}^s_{\tau'}]\ \ \mathbb{P}- a.s.$$
Now, using the fact that $\tilde{v}_3\in \mathcal{V}(s, \tau _1)$ provides a super-martingale inequality concatenated with any previous strategy $v$ and against any strategy $u$ of the opponent, we have
$$w^{\delta} (\tau _1', X_{\tau _1'})=w(\tau _1', X_{\tau _1'})\geq \mathbb{E}[w(\rho', X_{\rho'})|\mathcal{F}^s_{\tau _1'}] \geq \ \mathbb{E}[w^{\delta}(\rho', X_{\rho'})|\mathcal{F}^s_{\tau _1'}]\ \mathbb{P}- a.s.$$
The iterating conditioning provides 
$w^\delta (\tau ', X_{\tau '})\geq \mathbb{E}[w^\delta(\rho', X_{\rho'})|\mathcal{F}^s_{\tau'}]\ \ \mathbb{P}- a.s.$
This means $w^{\delta}\in \mathcal{U}^+$, so we  reached a contradiction.

1.2. The terminal condition property for $v^+$: the proof is quite similar to the proof of 1.1. More precisely, first we argue by contradiction, similar to either above or to \cite{bs-3} (but easier since controls are compact), then we construct a strategy $\tilde{v}$ as above. Finally,  we follow the arguments above based on It\^o's Lemma and conditioning, to finish the proof. We develop these ideas below.

 Assume that, for some $x_0\in \mathbb{R}^d$ we have 
$v^+(T,x_0)>g(x_0).$
 Since $g$ is USC , there exists an $\varepsilon >0$ such that
$g(x)\leq v^+(T,x_0)-\varepsilon$ if $ |(t,x)-(T,x_0)|\triangleq \max\{|x-x_0|,\  T-t \}\leq \varepsilon.$
Denote, once again, by 
$$\mathbb{T}\triangleq \overline{B(T,x_0,\varepsilon)}-B(T,x_0,\varepsilon /2))\cap( [0,T]\times \mathbb{R}^d).$$
Now $v^+$ is USC so bounded on $\mathbb{T}$.
Choose $\eta>0 $ small enough so that
\begin{equation*}
\label{bound1}
v^+(T,x_0)+\frac{\varepsilon^2}{4\eta} > \varepsilon +\sup_{(t,x)\in\mathbb{T}} v^+(t,x).
\end{equation*}
Since the inequality above is strict, we use again a decreasing sequence and a Dini argument from \cite{bs-2} to find some   $w\in \mathcal{V}^+$ such that
\begin{equation}
\label{bound}
v^+(T,x_0)+\frac{\varepsilon^2}{4\eta} > \varepsilon +\sup_{(t,x)\in \mathbb{T}} w(t,x).
\end{equation}
We now define, for $k>0$ the following function
$$\varphi^{\eta,\varepsilon, k}(t,x)=v^+(T,x_0)+\frac{|x-x_0|^2}{\eta}+ k(T-t). $$
Since $U$, $V$ are compact and $b, \sigma$ are continuous, 
for some  $k$ large enough (but no smaller than 
 $\varepsilon/2\eta$),  we have that 
$$-\varphi^{\varepsilon, \eta,k}_t- H^+(t,x, \varphi^{\varepsilon, \eta,k}_x, \varphi^{\varepsilon, \eta,k}_{xx})>0 \ \ \textrm{on}\ \  \overline{B(T,x_0,\varepsilon)}.$$
From \eqref{bound} we obtain
$\varphi^{\varepsilon, \eta,k}\geq \varepsilon + w$ on $\mathbb{T}$. 
Also, $$\varphi^{\varepsilon, \eta,k}(T,x)\geq v^+(T,x_0)\geq g(x)+\varepsilon \textrm{~for~} |x-x_0|\leq \varepsilon.$$  We  now  choose $\delta <\varepsilon$ and define 
$$v^{\varepsilon, \eta,k, \delta }=
\left \{
\begin{array}{l}
 w\wedge \Big ( \varphi ^{\varepsilon, \eta,k}-\delta  \Big)\textrm{~on~} \overline{ B(T, x_0, \varepsilon)},\\
w \textrm{~outside~}\overline{ B(T, x_0, \varepsilon)}.
\end{array}
\right.
$$
Using again the ideas in 
part 1.1  of the proof, we can show that $v^{\varepsilon, \eta, k, \delta}\in \mathcal{U}^+$ but $v^{\varepsilon, \eta, k, \delta}(T,x_0)=v^+(T,x_0)-\delta <v^+(T,x_0)$, leading to a contradiction. Therefore, $v^+(T, \cdot)\leq g$.

\noindent {\bf 2. \ \ $v^-$ is viscosity super-solution of upper Isaacs} 

2.  1. The interior super-solution property for $v^-$: Let $(t_0,x_0) \in [0,T)\times \mathbb{R}^d$ in the parabolic interior such that a smooth function $\varphi$ strictly touches $v^-$ from below at $(t_0,x_0)$. Assume, by contradiction, that
$\varphi _t +H^+(t,x,\varphi _x, \varphi _{xx}) >0$ at $(t_0, x_0)$. 
There exists a $\varepsilon >0$ and
there exists a measurable function $h:   V\rightarrow U,$
such that 
$$\varphi _t (t_0, x_0)+ \left[b(t_0,x_0,h(v), v)\cdot \varphi _x (t_0,x_0)+\frac{1}{2}Tr(\sigma(t_0,x_0,h(v),v)\sigma(t_0,x_0,h(v),{v})^T \varphi _{xx} (t_0,x_0)) \right] > \varepsilon.$$
In order to construct the function $h$ above, one does not need to use measurable selection arguments. Using uniform continuity on balls, one can actually choose a finitely valued $h$.
Since 
$\varphi _t  + b\cdot \varphi _x+\frac{1}{2}Tr(\sigma\sigma ^T \varphi _{xx}) $ is continuous in $(t,x,u,v)$ and $U,V$ are compact, it is therefore uniformly continuous  on 
$C\times U\times V$, where $C$ is any compact neighborhood of $(t_0,x_0)$.  
Therefore, we can choose a smaller $\varepsilon>0$ such that (despite the fact that $h$ is only measurable, and maybe not continuous) 
$$\varphi _t (t,x) + \left[b(t,x,h(v), v)\cdot \varphi _x(t,x)+\frac{1}{2}Tr(\sigma(t,x,h(v),v)\sigma(t,x,h(v),{v})^T \varphi _{xx}(t,x)) \right] > \varepsilon,\  {\textrm on}\ B(t_0, x_0, \varepsilon)\times V.$$
To be precise, small balls here are again intersected with the parabolic interior $[0,T)\times \mathbb{R}^d$.  
On the compact  torus $\mathbb{T}= \overline{B(t_0, x_0, \varepsilon)}- B(t_0, x_0, \varepsilon/2)$ we have that $\varphi <v^-$ and the max  of $\varphi -v^-$ is attained, therefore it is strictly negative. In other words
$\varphi +\eta <v^-$ on $\mathbb{T}$ for some $\eta >0$. Since  we can construct an increasing sequence of stochastic sub-solutions $w_n\nearrow v^-$, the same Dini type argument similar to \cite{bs-2} and \cite{bs-3} shows that, for $n$ large enough we have $\varphi +\eta/2 <w_n$ on $\mathbb{T}$. For simplicity, fix such an $n$ and call $w=w_n$. Now, define, for small $\delta <<\eta /2$
$$w^{\delta}\triangleq
\left \{ 
\begin{array}{ll}
(\varphi +\delta)\vee w\ \ {\textrm on}\ \  B(t_0, x_0, \varepsilon),\\
w \ \ {\textrm outside}\ \  B(t_0, x_0, \varepsilon).
\end{array}
\right.
$$
Since $w^{\delta}(t_0,x_0)>v^+(t_0,x_0)$, we obtain a contradiction if we can show $w^{\delta}\in \mathcal{U}^-$.  
In order to do so, fix $s$ and let $\tau:C([s,T])\rightarrow [s,T] $ be  a stopping time strategy for the  initial time $s$. Fix also $v\in \mathcal{V}(s,s)$. We need to construct an elementary  strategy $\tilde{u}\in \mathcal{U}(s,\tau)$ in the Definition of stochastic sub-solutions for $w^{\delta}$, depending on the $v$ we just fixed. This can be again done  as follows:
since $w$ is a stochastic sub-solution of upper Isaacs, there exists  an elementary strategy $\tilde{u}_1$ for $w$ starting at $\tau$ and depending on $v$ that does the job. Then
\begin{enumerate}
\item if $\varphi +\delta >w$ at $\tau$, follow the constant action $h(v).$
\item if $  w \geq \varphi +\delta $ at $\tau$ follow the strategy $\tilde{u}_1$
\item follow the  strategy  in 1-2 until the first time $\tau _1$ when $(t,X_t)\in \partial B(t_0, x_0, \varepsilon/2)$. On this boundary,  $w^{\delta}=w$
\item after this, follow the strategy $\tilde{u}_3\in \mathcal{U}(s, \tau _1)$ corresponding to the stochastic sub-solution $w$ with starting stopping rule $\tau _1$ and second player's strategy  $v\in \mathcal{V}(s,s).$
\end{enumerate}
In a rigorous notation though, we first define $\tilde{u}_2\in \mathcal{U}(s,\tau)$ by 
$$\tilde{u}_2(t,y (\cdot))= h (v(t, y(\cdot)) \, 1_{\{\varphi(\tau (y), y (\tau (y))) +\delta > w(\tau (y), y (\tau (y)))\}}+\tilde{u}_1 (t, y(\cdot))\, 1_{\{\varphi(\tau (y), y (\tau (y))) +\delta  \leq w(\tau (y), y (\tau (y)))\}}.$$
The key observation here, is that, since $h$ simply maps $V$ into $U$,  then, indeed, $\tilde{u}_2$ is an elementary strategy $\tilde{u}_2\in \mathcal{U}(s,\tau).$ Denote now by 
$\tau _1:C([s,T])\rightarrow [s,T]$, $\tau \leq \tau _1\leq T$, 
$$\tau _1(y)=\inf _{\tau (y)\leq t\leq T}\{(t,y(t))\in \partial B(t_0, x_0, \varepsilon/2)\}.$$
Now, if $\tilde{u}_3\in \mathcal{U}(s,\tau _1)$ is the strategy defining the stochastic super-solution $w$,  corresponding to $v\in \mathcal{V}(s,s)$ then we can define
$$\tilde{u}=\tilde{u}_2\otimes _{\tau _1}\tilde{u}_3\in \mathcal{U}(s, \tau)$$
to finish the proof. To actually check that $\tilde{u}$ (depending on $v$ and $\tau$) fulfills the condition in the definition of stochastic sub-solution for upper Isaacs, one has to follow similar arguments to the above part 1.1 of the proof, based on It\^o's Lemma  applied to $d(\varphi+\delta)(t, X_t)$, together with conditioning.

2.2. The terminal condition property for $v^-$: the proof is similar to the part 2.1 of the proof, but with an analytic construction very similar to the terminal condition part 1.2. More precisely, first we argue by contradiction, similar to the above step 1.2,   then we construct a strategy $\tilde{u}$ as in  2.1 above depending on the fixed $v$ and $\tau$. Finally,  we once again use  on It\^o's Lemma and conditioning in a similar manner to finish the proof.


\section{Verification by Comparison and  Dynamic Programming Principle}
\begin{Theorem}
Under the standing assumptions, we have
\begin{enumerate}
\item 
 $v^-=V^+=v^+$ is the unique continuous viscosity solution of the upper Isaacs equation. In addition, the upper value $V^+$ satisfies the (DPP) 
$$V^+(s,x)=\inf _{v\in \mathcal{V}(s,s)} \sup _{u\in \mathcal{U}(s,s)}\mathbb{E}\left [V^+(\rho(X^{s,x,u,v}_{\cdot}), X^{s,x,u,v}_{ \rho (X^{s,x,u,v}_{\cdot}   )} \right], \ \ \ \forall \rho \in \mathbb{B}^s.$$
\item  
$w^-=V^-=w^+$ is the unique continuous viscosity solution of the lower Isaacs equation. In addition, the lower value $V^-$  satisfies the (DPP) 
$$V^-(s,x)= \sup _{u\in \mathcal{U}(s,s)}\inf _{v\in \mathcal{V}(s,s)}  \mathbb{E}\left [V^-(\rho (X^{s,x,u,v}_{\cdot}), X^{s,x,u,v}_{ \rho (X^{s,x,u,v}_{\cdot}   )} \right], \ \ \ \forall \rho \in \mathbb{B}^s.$$

\end{enumerate}
If the Isaacs condition $H^-=H^+$ holds 
then
the game has a value and
 the value
 $$w^-=V^-=V=V^+=v^+$$
 is the unique continuous viscosity solution of the (unified) Isaacs equation. For each $\varepsilon >0$, there exist $\varepsilon$-saddle point strategies $(u,v)\in \mathcal{U}(s,s)\times \mathcal{V}(s,s)$.
\end{Theorem}

 We note that the DPP is  a statement about stopping  rules  $\rho\in \mathbb{B}^s$ and is {\emph not} formulated with respect   to stopping times $\rho '$ on the original probability space $(\Omega, \mathcal{F}, \mathbb{F}^s, \mathbb{P}).$  
\begin{Remark} In case the Isaacs condition holds, then one  only needs to prove half of the Theorem \ref{main1} (and arguably the easier one), namely that $v^+$ is a viscosity sub-solution of the Isaacs equation and that $w^-$ is a super-solution of the very same equation. These statements are perfectly symmetric (and this is the reason why we skipped completely the proof for $w^+$ and $w^-$). In this case, we already have that
$w^-\leq V^-\leq V^+-\leq v^+$
 and the comparison argument closes the proof (including the DPP).
\end{Remark}

Proof: the (DPP) is an easy conclusion of the way we defined $v^+, v^-, w^+$ and $w^-$, since the stochastic sub/super-solutions of upper/lower Isaacs satisfy the corresponding half DPP. In light of Theorem \ref{main1}, there is basically nothing left to prove, once a comparison result for bounded  semi-continuous viscosity solutions of the two Isaacs equations is available.  We make it available below to conclude $\diamond$

The  following Lemma  is basically borrowed from the analysis literature. One can either take it   from \cite{cil} (with modifications to allow for unbounded domains, bounded semi-continuous solutions and comparison up to $t=0$)  or (but still up to details) from  \cite{MR2533355}. 
\begin{Lemma}\label{lem:visc} Assume $b, \sigma$ are continuous, locally uniformly Lipschitz in $x$, satisfy linear growth (standing assumptions)  and $U$, $V$ are compact. Then
\begin{enumerate}
\item for each $R>0$ there exists some $\omega_R : [0, \infty)\rightarrow [0,\infty)$, with $\omega_R (0+)=0$ such that 
 for any 
$X,Y$ satisfying 
\begin{equation}\label{ishii}
-3 \alpha \left (\begin{array}{ll}I & 0\\0 &I
\end{array} \right )\leq  \left (\begin{array}{cc}X & 0\\0 &-Y
\end{array} \right )\leq 3\alpha  \left (\begin{array}{cc}I & -I\\-I&I
\end{array} \right )
\end{equation}
we have 
\begin{equation}\label{structural-lin}
L^{u,v}(t,y, \alpha (x-y),Y)-L^{u,v}(t,x,\alpha (x-y), X)\leq \omega_R (\alpha |x-y|^2+|x-y|), \ \forall \  t, u,v \ \textrm{and} \ |x|, |y|\leq R
\end{equation}
with the notation 
$L^{u,v}(t,x,p,M)\triangleq -b(t,x,u,v)\cdot p-\frac 12 Tr (\sigma (t,x,u,v)\sigma ^T(t,x,u,v)M)$
\item we have that, whenever condition \eqref{ishii} is satisfied,
\begin{equation}\label{structural-nonlinear}
F(t,y, \alpha (x-y),Y)-F(t,x,\alpha (x-y), X)\leq \omega _R(\alpha |x-y|^2+|x-y|), \ \forall \  t, u,v \ \textrm{and}\ |x|,|y|\leq R
\end{equation}
for $F=-H^+$ or $F=-H^-$
\item consider $v$ a bounded USC (on $[0,T]\times \mathbb{R}^d$) viscosity sub-solution of
$$-v_t+F(t,x,v_x, v_{xx})\leq 0\ \textrm{on} \ [0,T)\times \mathbb{R}^d$$
and $w$ a bounded  LSC (on $[0,T]\times \mathbb{R}^d$) viscosity super-solution of 
$$-w_t+F(t,x,w_x, w_{xx})\geq 0\ \textrm{on} \ [0,T)\times \mathbb{R}^d,$$
where (either) $F=-H^+$ (or $F=-H^-$).
If $v(T, \cdot)\leq w(T,\cdot)$ then $v\leq w$.

\end{enumerate}
\end{Lemma}
Proof: \begin{enumerate}
\item condition \eqref{structural-lin} (denoted by (3.14) in \cite{cil}) is actually proved in \cite{cil}  on page 19, in the elliptic case (no dependence on $t$), but the same proof applies to the parabolic case. One should have in mind that the local Lipschitz constant is uniform in $u,v$, therefore, the same  function $\omega_R$ can be chosen to work for all $u,v$
\item as still mentioned on page 19 in \cite{cil}, once \eqref{structural-lin} is uniform in $u,v$, it applies to any $\sup \inf$ or $\inf \sup$ transformation, obtaining the nonlinear structural condition \eqref{structural-nonlinear} (called, to recall, (3.14) in \cite{cil}) subject to \eqref{ishii}.
\item here, the proof is an adaptation of the comparison result in \cite{cil}  to   the case of an unbounded domain and (bounded)  semi-continuous solutions, under the structural condition \eqref{structural-nonlinear} subject to \eqref{ishii} and with the additional important piece of information that coefficients  have global linear growth. Comparison is obtained up to time $t=0$, as in \cite{MR2533355} or \cite{dfo}.   The case we have here does not fit perfectly in the existing literature, so that we cannot simply quote a single theorem, but the sequence of arguments used is well understood in the viscosity solution literature (see, for example, \cite{dfo} or \cite{MR2533355}. We  present the key arguments below.

To summarize,  we  use the penalization in 
the text-book \cite{MR2533355}, taking advantage of  the linear growth of the coefficients and the special structure of the HJB equation to reduce the problem to a bounded domain, and then appeal to a modification of the results in \cite{cil} based on  the arguments in  \cite{dfo} or \cite{MR2533355}, to obtain  comparison up to $t=0$ under the structural condition (3.14) in \cite{cil}.

Assume, for example, that $F=-H^+$. Let $\phi (t,x)=e^{-\lambda t} (1+|x|^2)$.
From the linear growth assumption {\bf{LG}} we have 
$$\sup _{u\in U}\sup _{v\in V}\{|b(t,x,u,v)|+|\sigma (t,x, u,v)|\}\leq C(1+|x|),\ \ \forall \ 0\leq t\leq T ,x\in \mathbb{R}^d.$$
This means that, for $\lambda >0$ large enough, we have
$$\phi _t(t,x)+\sup _{u\in U}\sup_{v\in V} [b(t,x,u,v)\cdot \phi _x(t,x)+\frac 12 Tr (\sigma (t,x,u,v)\sigma ^T(t,x,u,v)\phi _{xx}(t,x)) ]<0,$$
 for all $\ 0\leq t\leq T ,x\in \mathbb{R}^d.$ Therefore, for any $\varepsilon >0$, the LSC function 
 $w^{\varepsilon}=w+\varepsilon \phi$ is (easy to check) a viscosity super-solution of 
 $w^{\varepsilon}_t+H^+(t,x,w ^{\varepsilon}_x, w^{\varepsilon}_{xx})\leq 0.$
 According to  Theorem 8.2 in \cite{cil} (actually the follow up observations on page 52),  because $F=-H^+$  satisfies the structural condition \eqref{structural-nonlinear} subject   to \eqref{ishii}, we have, for any bounded open domain $\mathcal{O}=\{|x|<R\}\subset \mathbb{R}^d$,  a comparison result \emph{after time } $t=0$
 $$v(t,x)\leq w^{\varepsilon}(t,x)+\max \left ( \sup _{(0,T]\times \partial \mathcal{O}}(v-w^{\varepsilon})^+, \sup _{\{T\}\times \overline{\mathcal{O}}} (v-w^{\varepsilon})^+
 \right )\ \ \forall \ 0<t\leq T ,x\in \mathcal{O}.$$
However, it is well known (for example from \cite{dfo}, or \cite{MR2533355}), that, as long as the viscosity semi-solution property holds at time $t=0$ as well (which is the case here), the doubling argument  can be reproduced to provide comparison \emph{ up to time} $t=0$, i.e.
 $$v(t,x)\leq w^{\varepsilon}(t,x)+\max \left ( \sup _{[0,T]\times \partial \mathcal{O}}(v-w^{\varepsilon})^+, \sup _{\{T\}\times  \overline{\mathcal{O}}} (v-w^{\varepsilon})^+
 \right )\ \ \forall \ 0\leq t\leq T ,x\in \mathcal{O}.$$
  If one does not want to use the version of the Theorem of Sums (TOS) in \cite{dfo} to get comparison up to $t=0$ in the relation above, an alternative way is to simply define the whole game starting at an earlier negative time.
 Since both $v, w$ are bounded, the choice of $\phi$ ensures that, for large enough domain $\mathcal{O}=\{|x|<R\}$ we have $ \sup _{[0,T]\times \partial \mathcal{O}}(v-w^{\varepsilon})^+=0$.  In addition, since
 $v(T, \cdot)\leq w(T, \cdot)\leq w^{\varepsilon}(T, \cdot),$ we conclude that $\sup _{\{T\}\times  \overline{\mathcal{O}}} (v-w^{\varepsilon})^+=0$ as well.  In summary, for each $\varepsilon >0$ we have
 $$v(t,x)\leq w(t,x)+\varepsilon \phi (t,x), \ \forall \ 0\leq t\leq T, |x|< R$$
 for any $R$ large enough (depending on $\varepsilon$). We first let $R\nearrow \infty $ and then $\varepsilon \searrow 0$ to conclude that $v\leq w$ on $[0,T]\times \mathbb{R}^d$.

\end{enumerate}
\section{Conclusions}
Using Stochastic Perron's method  we  show that the values of differential games formulated symmetrically over elementary feed-back strategies are (the unique continuous) viscosity solutions of the corresponding Isaacs equations and the DPP holds. 
We  allow for the state equation to be degenerate and time-dependent.
Some of the assumptions (for example, compactness of $U$, $V$) seem to be removable, but one still needs enough to get comparison of semi-continuous viscosity solutions.

The program we propose amounts to verification for non-smooth viscosity solutions \emph{for the game in the original strong definition}, and  seems rather elementary compared to either the verification approach in \cite{swiech-2} or proving (parts of) the DPP working with the value functions as in \cite{fs}, \cite{MR1312585} or \cite{cr}.
We introduce the  class of Elementary (feed-back) Strategies, which have  strong intuition behind and are  easily amenable to analysis.  In particular, a player in a game that is believed to be continuous in time, is most likely to use Elementary Strategies to approximate the saddle point, if such exists. Obviously, one does not expect a saddle point to consist of elementary strategies, but, if the Isaacs condition holds, one has $\varepsilon$-saddle points among elementary strategies. The  present work sets the foundation  for the  analysis (in \cite{sirbu-2}) of more realistic zero-sum differential games over  mixed (elementary) feed-back strategies. Such a  game always has a value, which is  the unique solution of the (mixed) Isaacs equation.

We believe the present analysis sheds more light on the  connection between controlled diffusions and viscosity solutions of  fully non-linear PDE's. One important observation is that, in this symmetric formulation of the game over feed-back strategies,  the DPP is expected to hold only over stopping rules, i.e. stopping times where the decision to stop is based upon observing only the state, and not the noise. The technical parts of the proofs (which are basically elementary) are reminiscent of the previous proofs in \cite{bs-1}, \cite{bs-2} and especially \cite{bs-3}, and, for the analytic constructions, resemble the classic work of Ishii \cite{ishii}. An important part of the contribution resides  in the symmetric formulation of the game over elementary strategies (and strong solutions of the state equation) and in the proper (non-trivial) identification of stochastic semi-solutions to the Isaacs equations, that allows for the Stochastic Perron's  Method to be completed. With this in mind, we reiterate the conjecture that any stochastic optimization problem could be treated using Stochastic Perron's Method, provided that it   is properly formulated, and the stochastic semi-solutions defined accordingly.

\bibliographystyle{amsalpha} 
\providecommand{\bysame}{\leavevmode\hbox to3em{\hrulefill}\thinspace}
\providecommand{\MR}{\relax\ifhmode\unskip\space\fi MR }
\providecommand{\MRhref}[2]{%
  \href{http://www.ams.org/mathscinet-getitem?mr=#1}{#2}
}
\providecommand{\href}[2]{#2}

\end{document}